\documentclass{amsart}
\usepackage{amssymb,latexsym,subfigure}

\numberwithin{equation}{section}

\newtheorem{theorem}{Theorem}[section]
\newtheorem{proposition}[theorem]{Proposition}
\newtheorem{corollary}[theorem]{Corollary}
\newtheorem{remark}[theorem]{Remark}
\newtheorem{lemma}[theorem]{Lemma}

\newtheorem{definition}[theorem]{Definition}
\newtheorem{conjecture}[theorem]{Conjecture}

\def\proof{\smallskip\noindent {\bf Proof. }}
\def\endproof{\hfill$\square$\medskip}
\def\endproofmath{\quad\square}

\def\ii{\mathbf{i}}
\def\jj{\mathbf{j}}

\def\wnot{w_\mathrm{o}}
\def\RR{\mathbb{R}}
\def\ZZ{\mathbb{Z}}
\def\CC{\mathbb{C}}
\def\HH{\mathcal{H}}
\def\H{{\rm H}}
\def\l{\ell}

\def\Yo{Y^\circ}
\def\Yogu{Y_{\geq u}}

\def\str{\text{{\rm \texttt{str}}}}
\def\Card{\texttt{card}}

\def\lk{{\rm lk}}
\def\Lkuv{{\rm Lk}(u,v)}

\title[Stratified spaces formed by totally positive varieties]
{Stratified spaces formed by\\ totally positive varieties}

\author{Sergey Fomin}
\address{
Department of Mathematics, University of Michigan,
Ann Arbor, MI 48109, USA}
\email{fomin@math.lsa.umich.edu}

\author{Michael Shapiro}
\address{Department of Mathematics, 
Royal Institute of Technology, S-10044, Stockholm, Sweden} 
\email{mshapiro@math.kth.se} 

\thanks{The first author was supported in part 
by NSF grant \#DMS-9700927.
}

\date{December 14, 1999}


\begin{document}


\begin{abstract}
By a theorem of A.~Bj\"orner~\cite{bjorner-cw}, 
for every interval $[u,v]$ in the Bruhat order of a Coxeter group~$W$, 
there exists a stratified space
whose strata are labeled by the elements of $[u,v]$,  
adjacency is described by the Bruhat order, and 
each closed stratum (resp., the boundary of each stratum) 
has the homology of a ball (resp., of a sphere). 

Answering a question posed in~\cite{bjorner-cw}, 
we suggest a natural geometric realization of these stratified spaces  
for a Weyl group~$W$ of a semisimple Lie group~$G$, 
and prove its validity in the case of the symmetric group. 
Our stratified spaces arise as links in the 
Bruhat decomposition of the totally nonnegative part
of the unipotent radical of~$G$. 
\end{abstract}

\maketitle

\section{Introduction and main results}
\label{sec:intro}

In a 1984 paper~\cite{bjorner-cw}, A.~Bj\"orner has shown that
every  interval 
in the Bruhat order of a Coxeter
group~$W$ is the ``face poset'' of some stratified space, 
in which each closed stratum (resp., the boundary of each stratum)
has the homology of a ball (resp., of a sphere).
Passing to the Euler characteristic,  this result 
implies D.-N. Verma's formula~\cite{verma,verma-correction}
for the M\"obius function of the Bruhat order, viz.,  
$\mu(u,v)=(-1)^{\l(v)-\l(u)}$, $u\leq v$, 
where $\l$ denotes the length function.
(This is in turn equivalent to saying that each Bruhat interval 
contains equally many elements of even and odd length.)

Bj\"orner has in fact proved a stronger result,
namely, every interval in the Bruhat order is the face poset
of a regular cell complex
(i.e., closed strata actually \emph{are} balls). 
However, the construction of such complex 
in \cite{bjorner-cw} was entirely ``synthetic''
(essentially, a succession of cell attachments;
cf.\ \cite[4.7.23]{oriented-matroids}). 
Furthermore, it was based on the existence of a combinatorial
shelling, which by itself easily implies Verma's formula,
bypassing all geometry.
A question posed in~\cite{bjorner-cw} asked for
a natural geometric construction of 
a stratified space with the desired properties.

In this paper, we propose such a construction for the case where $W$ 
is the Weyl group of a semisimple group~$G$.  
In the type $A$ case, where $W$ is the symmetric group and $G$ the
special linear group,
we prove that our stratified spaces indeed have the required
homological properties. 
The spaces we construct are links of cells in the
Bruhat decomposition of the \emph{totally nonnegative} part
of the unipotent radical of~$G$.

In the remainder of Section~\ref{sec:intro},
we present the details of this construction,
and state our main results and conjectures.
The rest of the paper is devoted to proofs.


Let $G$ be 
a semisimple, simply connected
algebraic group defined and split over~$\RR$. 
Let $B$ and $B_-$ be two opposite Borel subgroups of $G$, so that
$H\!=\!B_-\cap B$  is an $\RR$-split maximal torus  in~$G$;
we denote by $N$ and $N_-$ the unipotent radicals of $B$ and $B_-\,$,
respectively.

For the type $A_{n-1}\,$,
the group $G$ is the real special linear group $SL(n,\RR)$;
$H$, $B$ and $B_-$ are the subgroups of diagonal,
upper-triangular, and lower-triangular matrices, respectively;
$N$ and $N_-$ are the subgroups of $B$ and $B_-$ that consist of
matrices whose diagonal entries are equal to~1.

We denote by $Y$ the set of all \emph{totally nonnegative} elements
in~$N$.
In the case of the special linear group,
$Y$ consists of the upper-triangular unipotent matrices
whose all minors are nonnegative.
The general definition was first suggested by G.~Lusztig
(see \cite{lusztig-survey} and references therein). 
In our current notation, Lusztig defined $Y$ 
as the multiplicative submonoid of $N$ generated by the elements
$\exp(te_i)$, $t\geq 0$,
where the $e_i$ are the Chevalley generators of the Lie
algebra of~$N$.
An alternative description in terms of nonnegativity of certain
``generalized minors'' was given in~\cite{FZosc} (cf.\ 
Proposition~\ref{prop:TP-via-gen-minors} below).

Let $W$ be the Weyl group of~$G$.
The length of an element $w\!\in\! W$ is denoted by~$\l(w)$. 
The group $W$ is partially ordered by the Bruhat order,
defined geometrically by 
$u\leq v\Longleftrightarrow B_- u B_-\subset \overline{B_- v B_-}\,$.
The Bruhat decomposition $G=\bigcup_{w\in W} B_- w B_-$ induces
the partition of $Y$ into mutually disjoint
\emph{totally positive varieties}
$\Yo_w = Y \cap B_- w B_-\,$, $w\in W$
(this terminology is borrowed from~\cite{FZ}).

We denote $Y_w=\overline{\Yo_w}$.
The varieties $\Yo_w$ 
were first studied by Lusztig
in~\cite{lusztig-reductive}, 
where, in particular, the following basic properties were obtained.

\begin{proposition}
\cite{lusztig-reductive}
\label{prop:lusztig-cells}
Each totally positive variety $\Yo_w$ is a cell;
more precisely, $\Yo_w$~is  homeomorphic to $\RR^{\l(w)}$.
Furthermore, $Y_w=\bigcup_{u\leq w}\Yo_u\,$.
\end{proposition}

\subsection*{Example: $G=SL(3,\RR)$}
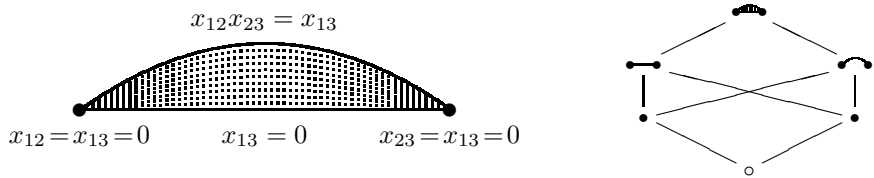
\begin{figure}[ht]
\centering
         \begin{picture}(140,45)(-10,-23)
\thicklines
\qbezier(0,0)(70,50)(140,0)
\qbezier[60](0,0)(70,45)(140,0)
\qbezier[60](0,0)(70,40)(140,0)
\qbezier[60](0,0)(70,35)(140,0)
\qbezier[60](0,0)(70,30)(140,0)
\qbezier[60](0,0)(70,25)(140,0)
\qbezier[60](0,0)(70,20)(140,0)
\qbezier[60](0,0)(70,15)(140,0)
\qbezier[60](0,0)(70,10)(140,0)
\qbezier[60](0,0)(70,5)(140,0)
\put(0,0){\line(1,0){140}}
\put(0,0){\circle*{5}}
\put(140,0){\circle*{5}}
\put(70,-10){\makebox(0,0){$x_{13}=0$}}
\put(70,33){\makebox(0,0){$x_{12}x_{23}=x_{13}$}}
\put(0,-10){\makebox(0,0){$x_{12}\!=\!x_{13}\!=\!0$}}
\put(140,-10){\makebox(0,0){$x_{23}\!=\!x_{13}\!=\!0$}}
         \end{picture}
\qquad\qquad\qquad\qquad
         \begin{picture}(80,65)(-40,0)
\put(0,0){\circle{3}}
\put(40,20){\circle*{3}}
\put(-40,20){\circle*{3}}
\put(-45,40){\circle*{3}}
\put(-45,40){\line(1,0){10}}
\put(-35,40){\circle*{3}}
\put(35,40){\circle*{3}}
\qbezier(35,40)(40,45)(45,40)
\put(45,40){\circle*{3}}
\put(-5,60){\circle*{3}}
\qbezier(-5,60)(0,65)(5,60)
\qbezier[10](-5,60)(0,64)(5,60)
\qbezier[10](-5,60)(0,63)(5,60)
\qbezier[10](-5,60)(0,62)(5,60)
\qbezier[10](-5,60)(0,61)(5,60)
\put(-5,60){\line(1,0){10}}
\put(5,60){\circle*{3}}
\thinlines
\put(-5,3){\line(-2,1){30}}
\put(5,3){\line(2,1){30}}
\put(-35,21){\line(4,1){65}}
\put(35,21){\line(-4,1){65}}
\put(-40,24){\line(0,1){12}}
\put(40,24){\line(0,1){12}}
\put(33,44){\line(-2,1){24}}
\put(-33,44){\line(2,1){24}}
         \end{picture}
\caption{Totally nonnegative varieties $\Yo_w$ 
in the special case $G=SL(3,\RR)$}
\label{fig:SL3}
\end{figure}

In this case,
\[
Y=\left\{
x=\left[
\begin{array}{ccc}
1 & x_{12} & x_{13} \\
0 & 1 & x_{23} \\
0 & 0 & 1
\end{array}
\right]
\,:\,
x_{12}\geq 0, ~ 
x_{23}\geq 0, ~ 
x_{13}\geq 0, ~ 
\left|\!\!
\begin{array}{ccc}
x_{12} & x_{13} \\
1 & x_{23} 
\end{array}
\!\right|
\geq 0
\right\}\,. 
\]
Thus the set $Y$ is described  
in the coordinates $(x_{12}, x_{23}, x_{13})$ 
as the closure of 
one of the pieces into which the plane $x_{13}\!=\!0$ 
and the hyperbolic paraboloid
$x_{12}x_{23}\!=\!x_{13}$ partition the 3-space,---namely, 
the piece containing the point $(1,1,\frac{1}{2})$. 
The semialgebraic set~$Y$ decomposes naturally into 6 algebraic strata:
the origin, two~rays (the positive semi-axes for $x_{12}$ and
$x_{23}$), two 2-dimensional pieces connecting them,
and the 3-dimensional interior. 
These are the 6 Bruhat strata $\Yo_w\,$, 
for $w\in W=\mathcal{S}_3$ (the symmetric group). 
Figure~\ref{fig:SL3} shows a planar cross-section of this 
stratification---or, equivalently, the link of the $0$-dimensional cell. 
The adjacency of the strata $\Yo_w$
is indeed described by the Bruhat order on
$\mathcal{S}_3$, in agreement with
Proposition~\ref{prop:lusztig-cells}.

For $u,v\in W$, $u\leq v$,
the \emph{Bruhat interval} $[u,v]$ is defined by
$[u,v]=\{u\leq w\leq v\}$, with the partial order inherited from~$W$.
Similarly, $(u,v]\stackrel{\rm def}{=}\{u< w\leq v\}$. 

In view of Proposition~\ref{prop:lusztig-cells}, it is natural to suggest 
that the geometric model for a Bruhat interval $[u,v]$ 
(or $(u,v]$) 
is provided by the \emph{link} 
\[
\Lkuv\!=\!\lk(\Yo_u,Y_v) 
\]
of the cell $\Yo_u$ inside the subcomplex $Y_v\subset Y$. 

The following are our main results. 

\begin{theorem}
\label{th:strata}
For any $u\leq v$, the link $\Lkuv$ is well defined as a
stratified space.
The strata $S_{u,v,w}=\Lkuv\cap \Yo_w$ are labelled by the elements
$w\in (u,v]$, and  
each  stratum $S_{u,v,w}$ is an open smooth manifold
of dimension $\l(w)-\l(u)-1$. 
The closures and boundaries of the strata $S_{u,v,w}$ are given by
\begin{equation}
\label{eq:closure+boundary}
\overline{S_{u,v,w}}=\displaystyle\bigcup_{u< w'\leq w}
S_{u,v,w'}\,,
\qquad
\partial \overline{S_{u,v,w}} 
= \displaystyle\bigcup_{u< w'<w} S_{u,v,w'}\,.
\end{equation}
\end{theorem}

(Here by a stratified space we mean a decomposition 
of a semialgebraic set into disjoint smooth submanifolds
labelled by the elements of a partially ordered set,
as described for example in \cite[Section~1.1]{GM}. 
Although the stratifications we consider seem to satisfy Whitney's regularity
conditions (cf.\ \cite[Section~1.2]{GM}), we will not need to verify
these conditions to justify our constructions.)

\begin{theorem} 
\label{th:orientable}
{\rm (Type~$A$ only.)}
All the strata $S_{u,v,w}$ are orientable.
\end{theorem}

\begin{theorem}
\label{th:contractions}
{\rm (Type~$A$ only.)}
Each closed stratum $\overline{S_{u,v,w}}$ is contractible.
Moreover, the contraction can be chosen so that 
it restricted to a contraction of $S_{u,v,w}\,$. 
\end{theorem}

These theorems ensure that the stratified spaces $\Lkuv$
have the desired homological properties, as we will now explain. 
Let $\H_i(X)$ (resp., $\H_i(X,X')$) denote, as usual, 
the ordinary $i$-th homology group 
of a CW-complex~$X$ 
(resp., pair of CW-complexes  $X'\subset X$). 
The corresponding Euler characteristics are denoted by 
$\chi(X)$ and $\chi(X,X')$, respectively. 

\begin{corollary}
\label{cor:H(S,dS)}
{\rm (Type~$A$ only.)}
For any $u<w\leq v$, we have 
\begin{equation}
\label{eq:H(S,dS)}
\H_i(\overline{S_{u,v,w}},\partial \overline{S_{u,v,w}})=
\left\{
\begin{array}{l}
\mathbb{Z},\,\text{ if } i=\l(w)-\l(u)-1\\
0,\,\text{ otherwise. }\\
\end{array}
\right.
\end{equation}
Consequently, $\chi(\overline{S_{u,v,w}},\partial
\overline{S_{u,v,w}})=(-1)^{\l(w)-\l(u)-1}$. 
\end{corollary}

\begin{proof}
We will need the Lefschetz duality isomorphism
\cite[Exercise~18.3]{DFN}:
\[
\H_i(X,A)\simeq \H^{n-i}(X\setminus A;\ZZ)\,, \qquad i>0\,, 
\]
where $X$ is a compact topological space
and $A$ its closed subset such
that $X\setminus A$ is a smooth orientable
$n$-dimensional manifold. 
Take $X=\overline{S_{u,v,w}}$ and 
$A=\partial\overline{S_{u,v,w}}$.
Then Theorems~\ref{th:strata} and~\ref{th:orientable} 
ensure that the above conditions are satisfied,
with $n=\l(w)-\l(u)-1$. 
Hence $\H_i(\overline{S_{u,v,w}},\partial \overline{S_{u,v,w}})=
\H^{\l(w)-\l(u)-1-i}( S_{u,v,w})$. 
Since 
$S_{u,v,w}$ is contractible by Theorem~\ref{th:contractions}, 
(\ref{eq:H(S,dS)}) follows. 
\end{proof}

\begin{corollary}
\label{cor:main}
{\rm (Verma's Theorem)}
\ $\displaystyle\sum_{u\leq w\leq v} (-1)^{l(w)}=0$.
\end{corollary}

Thus every Bruhat interval is an \emph{Eulerian poset}~\cite{StaNato}.

\begin{proof}
The additivity of the Euler characteristic \cite[V.5.7]{dold}, 
which applies in view of Theorem~\ref{th:strata}, gives 
\[
\chi(\Lkuv)= \sum_{u < w\leq v}
\chi(\overline{S_{u,v,w}},\partial \overline{S_{u,v,w}}).
\]
In the last identity, the left-hand side is equal to~1 
by Theorem~\ref{th:contractions}, 
while the right-hand side is equal to 
$\sum_{u < w\leq v} (-1)^{l(w)-l(u)-1}$ 
by Corollary~\ref{cor:H(S,dS)}. 
Simplifying, we obtain the desired formula. 
\end{proof}


For the type $A$, we prove the following refinement of
Theorem~\ref{th:strata}.
Let us define the stratified space $Y_{[u,v]}$ by
$Y_{[u,v]}=\bigcup_{w\in [u,v]}\Yo_w\,$.
Note that $\Lkuv\!=\!\lk(\Yo_u,Y_v)\!=\!\lk(\Yo_u,Y_{[u,v]})$. 

\begin{theorem}
\label{th:product}
{\rm (Type~$A$ only.)}
The stratified space $Y_{[u,v]}$ has the structure of the direct product
of the cell $\Yo_u$ and the cone over the link~$\Lkuv$.
More precisely, there exists an isomorphism of stratified spaces
$Y_{[u,v]}$ and $\Yo_u\times {\rm Cone}(\Lkuv)$,
whose restriction to each stratum is a diffeomorphism. 
\end{theorem}

\begin{remark}
{\rm 
Theorem~\ref{th:orientable} can be deduced from
Theorem~\ref{th:product}, as follows. 
By (\ref{eq:closure+boundary}), the stratum $S_{u,v,w}$ coincides with
the interior of ${\rm Lk}(u,w)$. 
Thus Theorem~\ref{th:product} asserts, in particular, that the cell
$\Yo_w$ is a direct product of the cell $\Yo_u$ and the interior of
the cone over ${\rm Lk}(u,w)$. Both cells $\Yo_w$ and $\Yo_u$ are
evidently orientable. Therefore (see, e.g.,
\cite[Exercise~3.2.24]{GP}) the interior of the cone over
${\rm Lk}(u,w)$ is orientable and so is the interior of
${\rm Lk}(u,w)$.
\qed
}
\end{remark}


\begin{conjecture}
\label{conj:general}
Theorems~\ref{th:contractions} and~\ref{th:product} 
(hence Theorem~\ref{th:orientable} and Corollary~\ref{cor:H(S,dS)}) 
hold 
for any semisimple algebraic group~$G$.
\end{conjecture}

We believe that Conjecture~\ref{conj:general} can be strengthened as
follows.

\begin{conjecture}
\label{conj:main}
Each stratum $S_{u,v,w}$
(resp., its closure, its boundary)
is homeomorphic to an affine space
(resp., closed ball, a sphere)
of dimension $\l(v)-\l(u)-1$
(resp., $\l(v)-\l(u)-1$, $\l(v)-\l(u)-2$).
Thus $\Lkuv$ is a regular cell complex.
\end{conjecture}

Assuming Conjecture~\ref{conj:main} holds,
each stratified link $\Lkuv$ provides a geometric realization
of the ``generalized synthetic Schubert variety'' whose existence was
hypothesized by Bj\"orner~\cite{bjorner-cw}.


We hope to extend the construction of the spaces $\Lkuv$ 
to an arbitrary simply-laced Coxeter group, and possibly further, 
so that the analogues of all statements 
formulated above would still hold. 
(Note that Bj\"orner's original result applies to intervals
in \emph{any} Coxeter group.) 

It should be mentioned that one of our ``hidden motivations'' has been
the desire to better understand the combinatorics of
Kazhdan-Lusztig polynomials.
It was already pointed out in the original paper by Kazhdan and
Lusztig~\cite{KL} that Verma's formula is equivalent to the assertion
that the constant term of any Kazhdan-Lusztig polynomial is~1.

\medskip

The remainder of this paper is organized as follows. 
Sections~\ref{sec:prelim}--\ref{sec:projection} introduce some 
useful Lie-theoretic machinery;
in particular, we define a projection onto a cell~$\Yo_u$ that
plays a crucial role in subsequent proofs. 
In Section~\ref{sec:links}, we prove 
Theorem~\ref{th:strata}. 
Section~\ref{sec:proofs-type-A} contains the proofs 
of Theorems~\ref{th:contractions} and~\ref{th:product}. 
These proofs are based on a technical lemma
(Lemma~\ref{lem:link}), which is proved in 
Section~\ref{sec:main-lemma} for the special case of $G=SL(n)$;
this is the only ``type-specific'' ingredient of our proofs. 

\medskip

\textsc{Acknowledgments.}
The authors are grateful to Saugata Basu, Ilia Itenberg, Viatcheslav
Kharlamov, Boris Shapiro, Eugenii Shustin, Viktor Vassiliev, and
Andrei Zelevinsky for valuable advice.  

\section{Preliminaries}
\label{sec:prelim}

This section introduces necessary technical background;
throughout it, we do not claim any originality. 

\nobreak

The notation used below is consistent with~\cite{FZ}.
In particular, we denote~by 
\[
G_0=B_-B=N_-HN
\] 
the set of elements of $x\in G$ that have a
Gaussian decomposition;
for the latter, we use the notation $x=[x]_-[x]_0[x]_+\,$. 

We think of the Weyl group $W$ as the quotient of the normalizer of
$H$ modulo~$H$, 
and identify each element $w\in W$ with a fixed 
representative in~$G$. 

\begin{lemma}
\label{lem:wBw-in-G0}
For $w\in W$, 
we have $w^{-1}B_- w\subset G_0$ and $w^{-1}Bw\subset G_0$ .
Moreover, $w^{-1}N_- w\subset N_- N$ and $w^{-1}Nw\subset N_- N$. 
\end{lemma}

\begin{proof}
Since $B = H N$, $B_- = H N_-\,$, and $w$ normalizes~$H$, 
it suffices to prove the last statement. 
It is well known (cf., e.g., \cite[Proposition~2.12]{FZ}
or \cite[(5.3)]{BZ})
that any $x\in N$ is uniquely factored as $x=x_1 x_2$
with $x_1\in N \cap w N_- w^{-1}$ and $x_2\in  N \cap w N w^{-1}$. 
Hence $N\subset w N_- w^{-1} \cdot w N w^{-1} = w N_-Nw^{-1}$, 
as desired. 
\end{proof}

\begin{lemma}
\label{lem:[]+insubgroup}
If $z\in w^{-1} B_- w$,
then $[z]_-\,,[z]_+\in w^{-1}N_- w$. 
Analogously, if $z\in w^{-1} B w$,
then $[z]_-\,,[z]_+\in w^{-1}N w$. 
\end{lemma}

\begin{proof}
It is enough to show that $z\in w^{-1} B w$
implies $[z]_+\in w^{-1}N w$. 
Just as in the proof of Lemma~\ref{lem:wBw-in-G0},
we can write $z=h w^{-1} x_1 x_2 w$,
where $h\in H$, $x_1\in N \cap w N_- w^{-1}$, 
and $x_2\in  N \cap w N w^{-1}$. 
Then $z=(h w^{-1}x_1 w)(w^{-1}x_2 w)$,
where the factors belongs to $B_-$ and $N$,
respectively. 
Thus $[z]_+=w^{-1}x_2 w\in w^{-1}N w$, as desired. 
\end{proof}

We define the subgroups 
\[
\begin{array}{rll}
N_-(w)&=w^{-1}B w\cap N_-&=w^{-1}N w\cap N_- \,,\\[.1in]
  N(w)&=w^{-1}B w\cap N  &=w^{-1}N w\cap N 
\end{array}
\]
and the set 
\[
N^w  = B_-wB_-\cap N \,.
\]

\begin{lemma}
\label{lem:x and y}
For any $w\in W$ and $x_w\in N^w$, 
there exists a unique $y\in N_-(w)$ satisfying $x_w=[wy]_+\,$.
Specifically, $y=w^{-1}[x_w w^{-1}]_+w\,$. 
\end{lemma}

\begin{proof}
Immediate from \cite[Propositions 2.10 and 2.17]{FZ}. 
\end{proof}

\begin{lemma}
\label{lem:action}
Let $x_w\in N^w$, $b_-\in B_-\,$, and $x_w b_-\in G_0\,$.
Then $[x_w b_-]_+\in N^w$. 
\end{lemma}

\begin{proof}
$[x_w b_-]_+\in B_- x_w b_- \subset B_- \cdot B_-wB_- \cdot b_-
=B_-wB_-\,$. 
\end{proof}

The following statement, however obvious, is quite useful. 

\begin{lemma}
\label{lem:trivial}
If $x\in G_0$ and $y\in G$, then $[[x]_+y]_+=[xy]_+$,
provided one of the two sides is well defined. \hfill \qed
\end{lemma}

\begin{lemma}
\label{lem:recovering-shift}
For any $x_w\,,\tilde x_w\in N^w\,$, 
there exists a unique $n_1\in N_-(w)$ satisfying $x_w=[\tilde x_w
n_1]_+\,$.  
Specifically, 
$n_1=w^{-1}([\tilde x_w w^{-1}]_+)^{-1}[x_w w^{-1}]_+w$. 
\end{lemma}

\begin{proof}
Uniqueness follows from the uniqueness part of 
Lemma~\ref{lem:x and y}, together with Lemma~\ref{lem:trivial}
and the fact that $N_-(w)$ is a group.
In more detail: assume $x_w=[\tilde x_w n_1]_+=[\tilde x_w n'_1]_+\,$,
where $n_1\neq n'_1$ and $n_1,n'_1\in N_-(w)$.
Let $y$ be as in Lemma~\ref{lem:x and y}. 
Then 
 $x_w=[\tilde x_w n_1\cdot n_1^{-1}n'_1]_+
=[x_w n_1^{-1}n'_1]_+
=[wy n_1^{-1}n'_1]_+\,$,
where $y\neq y n_1^{-1}n'_1\in N_-(w)$, a contradiction. 

With the notation $y=w^{-1}[x_w w^{-1}]_+w$
and $\tilde y=w^{-1}[\tilde x_w w^{-1}]_+w$,
it remains to check that $n_1\!=\!\tilde y^{-1} y$
satisfies $x_w\!=\![\tilde x_w n_1]_+\,$.  
Indeed, $x_w=[wy]_+=[w\tilde y n_1]_+=[x_w n_1]_+\,$. 
\end{proof}

We now turn to total nonnegativity.
Let us first recall Lusztig's original
definition~\cite{lusztig-reductive}. 
According to it, the set $Y$ of totally nonnegative elements in $N$
is defined as 
the multiplicative monoid generated by the elements
\begin{equation}
\label{eq:x_i(t)}
x_i(t)=\exp(te_i),
\end{equation}
where $t\geq 0$ 
and the $e_i$ are the Chevalley generators of the Lie algebra of~$N$.
One of the first results in~\cite{lusztig-reductive} is the following
description of the
Bruhat stratum $\Yo_w=Y\cap B_- w B_-\,$.

\begin{proposition}
\label{prop:param-stratum}
\cite{lusztig-reductive}
Let $(a_1,\dots,a_l)$ be a reduced word for~$w\in W$. 
Then the map 
\begin{equation}
\label{eq:param-stratum}
(t_1,\dots,t_l) \mapsto x_{a_1}(t_1)\cdots x_{a_l}(t_l)
\end{equation}
is a bijection between $\RR_{>0}^l$ and $\Yo_w\,$.  
\end{proposition}

(It is clear from this description that $\Yo_w$ is indeed a
cell of dimension~$l=\l(w)$; cf.\ Proposition~\ref{prop:lusztig-cells}.)
One is tempted to use the parametrizations
(\ref{eq:param-stratum}) to prove our main theorems.
Unfortunately, this approach encounters
substantial difficulties, chiefly due to the fact that 
the relationship between parametrizations of
adjacent cells is generally quite complicated. 
In what follows, we hardly make any use of
Proposition~\ref{prop:param-stratum}.

\begin{proposition}
{\rm (G.~Lusztig~\cite[6.3]{lusztig-survey})}
\label{prop:connected component}
The cell $\Yo_w$ is a connected component in~$N^w$
(in the ordinary topology).
\end{proposition} 


For $u\in W$, we denote $\Yogu=\bigcup_{v\geq u} \Yo_v$. 

To state the next result, we will need the notion of a
generalized minor of an element $x\in G$,
for which the reader is referred to~\cite[Section~1.4]{FZ}. 
Generalized minors are certain regular functions on~$G$, which can be
defined as suitably normalized matrix coefficients corresponding 
to pairs of extremal weights in some fundamental
representation of~$G$.
In the case of type~$A$, this notion coincides with the ordinary
notion of a minor of a square matrix. 

\begin{proposition}\cite[Theorem~3.1]{FZosc}
\label{prop:TP-via-gen-minors}
An element $x\in G$ is totally nonnegative, in the sense of
Lusztig~\cite{lusztig-reductive}, if and only if all its generalized
minors are nonnegative.  
\end{proposition}

\begin{lemma}
\label{lem:genericity}
If a generalized minor does not vanish at some point $x\in\Yo_u\,$, 
then it vanishes nowhere in $\Yo_u\,$,
and furthermore nowhere in $\Yogu$. 
\end{lemma}

\begin{proof}
For the type $A$, this is an immediate corollary
of \cite[Proposition~5.2.2]{BFZ}.
The general case can be deduced from (highly nontrivial) 
\cite[Proposition~7.4]{BZ-tensor-preprint}. 
According to the latter, 
for any generalized minor $\Delta$ and any sequence of indices
$a=(a_1,\dots,a_m)$, the function
$P_a(t_1,\dots,t_m)=\Delta(x_{a_1}(t_1)\cdots x_{a_m}(t_m))$ 
(cf.\ (\ref{eq:param-stratum})) is either identically zero,
or is a polynomial with positive integer
coefficients. 
(The type~$A$ version of this statement is well known;
see, e.g., \cite[Theorem~2.4.4]{BFZ}.) 
Since $\Delta$ does not vanish at some point in $\Yo_u\,$,
we know that $P_a$ is a nonzero polynomial for any reduced word $a$
for~$u$. 
For $v\geq u$, any reduced word $b$ for $v$ contains some reduced
word $a$ for $u$ as a subword (see \cite[5.10]{humphreys}). 
Hence $P_a$ is a specialization of $P_b\,$, obtained by setting some
of the variables equal to~$0$. 
Then $P_a\neq 0$ implies $P_b\neq 0$.
On the other hand, $P_b$ is a polynomial with positive coefficients,
so $P_b(t_1,t_2,\dots)\neq 0$ for any $t_1,t_2,\ldots >0$,
or. equivalently, $\Delta(x)\neq 0$ for any $x\in\Yo_v\,$. 
\end{proof}

\begin{lemma}
\label{lem:G0membership1}
For any $u\in W$, we have $B_-uB_- \subset G_0 u$. 
In particular, $N^u\subset G_0 u$. 
\end{lemma}

\begin{proof}
Follows from Lemma~\ref{lem:wBw-in-G0}. 
\end{proof}


\begin{corollary}
\label{cor:genericity-special}
$\Yogu \subset G_0 u$. 
\end{corollary}

\begin{proof}
By \cite[Corollary~2.5]{FZ}, the set 
$G_0 u$ is defined by several inequalities of the form
$\Delta\neq 0$, where $\Delta$ is a generalized minor. 
Since $\Yo_u\subset G_0 u$ (by Lemma~\ref{lem:G0membership1}),
none of these minors vanishes on $\Yo_u$---and therefore
none vanishes anywhere on $\Yogu\,$, by Lemma~\ref{lem:genericity}. 
\end{proof}

\begin{theorem}
{\rm (V.~V.~Deodhar~\cite[Corollary~1.2]{deodhar};
cf.~also \cite[Sec.~1.2]{KL2})}
\label{th:deodhar}
For $u,v\in W$, the intersection $B_-vB_-\cap BuB_-$ is non-empty if
and only if $u\leq v$. 
\end{theorem}

\begin{corollary}
\label{cor:deodhar}
$G_0 u\subset \bigcup_{v\geq u} B_-vB_-\,$. 
\end{corollary}

\begin{proof}
Let $x\in G_0 u$ and $x\in B_-vB_-\,$, $v\in W$.
Then, by Theorem~\ref{th:deodhar}, 
\[
B_-vB_-\cap B_-Bu\neq\emptyset 
\Longrightarrow B_-vB_-\cap Bu\neq\emptyset  
\Longrightarrow u\leq v \, ,
\]
as desired. 
\end{proof}

\begin{corollary}
\label{cor:Yu=Y-and-G0u}
$\Yogu=Y\cap G_0 u$. 
\end{corollary}

\begin{proof}
The inclusion $\Yogu\subset Y\cap G_0 u$ is
Corollary~\ref{cor:genericity-special}. 
The opposite inclusion is immediate from Corollary~\ref{cor:deodhar}. 
\end{proof}

\begin{lemma}
\label{lem:YY-in-Y}
$\Yogu \,Y\subset\Yogu\,$.
\end{lemma}

\begin{proof}
By \cite[Lemma~2.14]{lusztig-reductive},
for any $w_1,w_2\in W$, we have 
$\Yo_{w_1}\Yo_{w_2}=\Yo_{w_3}$ for some $w_3\in W$.
Moreover, it is clear from the proof of this statement
in~\cite{lusztig-reductive} that $w_3\geq w_1\,$,
and the lemma follows. 
\end{proof}

\subsection*{Example: $G=SL(3,\RR)$}
Let $u=s_1$, the transposition of 1 and 2 in the symmetric
group $W=\mathcal{S}_3\,$. 
Then,
using the notation 
$x=\left[
\begin{array}{ccc}
1 & x_{12} & x_{13} \\
0 & 1 & x_{23} \\
0 & 0 & 1
\end{array}
\right]$ for the elements $x\in N$, we have 
\[
N(u)=\{x_{12}=0\} 
\,, 
\quad 
N^u=\left\{
\begin{array}{cc}
x_{12}\neq 0 & x_{13}=0\\
         & x_{23}=0 
\end{array}
\right\}\,,\quad 
N\cap G_0u=
\{x_{12}\neq 0 \}  
\]
and
\[
\Yo_u=\left\{
\begin{array}{cc}
x_{12}>0 & x_{13}=0\\
         & x_{23}=0 
\end{array}
\right\}
\,,\quad 
\Yogu=
\left\{
\begin{array}{cc}
x_{12}>0 & x_{13}\geq 0\\
         & x_{23}\geq 0 
\end{array}
\ \ 
\left|\!\!
\begin{array}{ccc}
x_{12} & x_{13} \\
1 & x_{23} 
\end{array}
\!\right|
\geq 0
\right\}\,. 
\]

\section{Projecting on a 
cell}
\label{sec:projection}

In this section, we introduce a projection 
$\pi_u:\Yogu\to\Yo_u$ that will later be used to construct and study
the links $\Lkuv\!=\!\lk(\Yo_u,Y_v)$. 
This projection can be viewed as the totally positive version of
the projection of an affine open
neighborhood of a Schubert cell onto the cell itself, 
which arises from the direct product decomposition 
described by Kazhdan and Lusztig in \cite[Secs.~1.3--1.4]{KL2}

Let us fix an element $u\in W$. 

\begin{lemma}
\label{lem:G0membership2}
If $x\in G_0 u\cap G_0$
(in particular, if $x\in \Yogu$---cf.\
Corollary~\ref{cor:genericity-special}),  
then $u^{-1}[xu^{-1}]_+u\in G_0$ and 
$u[u^{-1}[xu^{-1}]_+u]_-\in G_0$. 
\end{lemma}

\begin{proof}
The first statement follows from Lemma~\ref{lem:wBw-in-G0}.
Proof of the second one: for some $b_-\!\in\! B_-$ and $b\!\in\! B$, we have 
$u[u^{-1}[xu^{-1}]_+u]_-
= u u^{-1}b_-xu^{-1}u b
=b_-xb\in G_0$. 
\end{proof}

\begin{lemma}
\label{lem:pi_u}
The map $(x_u, x^u)\mapsto x=x_u x^u$ is a bijection 
\[
N^u \times  N(u)\to N\cap G_0 u.
\]
The inverse map $x\mapsto (x_u,x^u)$ is given by 
\begin{equation}
\label{eq:x_u-formula}
x_u=[u[u^{-1}[xu^{-1}]_+u]_-]_+
\end{equation}
and 
\begin{equation}
\label{eq:x^u-formula}
x^u=[u^{-1}[x u^{-1}]_+u]_+ \,. 
\end{equation}
Furthermore, if $x\in N\cap G_0 u$ is totally nonnegative
(i.e., $x\in \Yogu$; cf.\ Corollary~\ref{cor:Yu=Y-and-G0u}),
then $x_u$ is totally nonnegative
(i.e., $x_u\in \Yo_u$). 
\end{lemma}

\begin{proof}
Assume $x_u\in N^u$, $x^u\in N(u)$, and $x=x_u x^u$.
Then 
$x=x_u x^u 
\in G_0 u\cdot u^{-1}N u = G_0 u$
(by Lemma~\ref{lem:G0membership1} and the definition of $N(u)$),
as claimed. 

Let us prove that the map in question is a surjection.
Let $x\in N\cap G_0u$, and let $x_u$ and $x^u$ be given by
(\ref{eq:x_u-formula}) and (\ref{eq:x^u-formula});  
note that the right-hand sides of these formulas 
are well defined (by Lemma~\ref{lem:G0membership2}). 
Thus $x_u=[u[y]_-]_+$ and $x^u=[y]_+\,$,
where $y=u^{-1}[xu^{-1}]_+u$. 
Then $x_u\in B_-u[y]_-\subset B_-uB_-$
and $x^u\in N(u)$ (by Lemma~\ref{lem:[]+insubgroup}).
Furthermore, $x_u x^u= [u[y]_-]_+[y]_+
=[uy]_+$ (since $y\in N_-N$ by Lemma~\ref{lem:wBw-in-G0})
and therefore $x_u x^u=[uy]_+=[[xu^{-1}]_+u]_+=x$
(by Lemma~\ref{lem:trivial}). 

Let us now prove injectivity. 
Suppose $x_u\in N^u$, $x^u\in N(u)$, and $x=x_u x^u$. 
We will show that $x_u$ and $x^u$ can be 
recovered from $x$ via 
(\ref{eq:x_u-formula})--(\ref{eq:x^u-formula}). 
Since $x\in N\cap G_0 u$, the right-hand sides of 
(\ref{eq:x_u-formula})--(\ref{eq:x^u-formula})
are well defined (by Lemma~\ref{lem:G0membership2}). 
Then
\[
\begin{array}{rcll}
[u[u^{-1}[xu^{-1}]_+u]_-]_+
&=& [u[u^{-1}[x_ux^uu^{-1}]_+u]_-]_+\\[.1in]
&=& [u[u^{-1}[x_u u^{-1}]_+ ux^u u^{-1} u]_-]_+ & 
    \textrm{(since $ux^u u^{-1}\in N$)}\\[.1in]
&=& [u[u^{-1}[x_u u^{-1}]_+ u]_-]_+\\[.1in]
&=& [u\cdot u^{-1}[x_u u^{-1}]_+ u]_+ & 
    \textrm{(by Lemma~\ref{lem:x and y})}\\[.1in]
&=& x_u & 
    \textrm{(by Lemma~\ref{lem:trivial}),}
\end{array}
\]
proving (\ref{eq:x_u-formula}).

Let us prove (\ref{eq:x^u-formula}).
Denote $A=u^{-1} [x u^{-1}]_+u$. 
We have 
\[
x=[x]_+=[u (u^{-1} [x u^{-1}]_+ u)]_+=[uA]_+
=[u[A]_-]_+[A]_+
\]
(by Lemma~\ref{lem:wBw-in-G0}). 
On the other hand, we already proved that $x_u=[u [A]_-]_+\,$.
Thus $x=x_u [A]_+$, i.e., $x^u=[A]_+\,$, as desired. 

It remains to prove that $x_u$ is totally nonnegative whenever $x$
is. Assume $x\in \Yo_w\subset\Yogu\,$. 
Consider a path that connects $x$ with a point $x_0\in \Yo_u$
and stays inside $\Yo_w$ (such a path exists since $\Yo_w$ is
connected and its boundary contains $\Yo_u$---see
Proposition~\ref{prop:lusztig-cells}). 
The image of this path under the projection 
$N\cap G_0 u\to N^u$ 
connects $x_u$ with $x_0\,$.
Since $x_0\in \Yo_u\,$, Proposition~\ref{prop:connected component} 
implies that $x_u\in \Yo_u\,$. 
\end{proof}

In view of Lemma~\ref{lem:pi_u}, the formula
\begin{equation}
\label{eq:pi_u}
\pi_u(x)=[u[u^{-1}[xu^{-1}]_+u]_-]_+   
\end{equation}
defines a continuous projection 
$\pi_u:\Yogu\to\Yo_u\,$. 
%
%
(The map $\pi_u$ is a projection 
since $x=x\cdot 1$ gives the factorization in question for $x\in\Yo_u\,$.)

\subsection*{Example: $G=SL(3,\RR)$, $u=s_1\,$} 
For $x\in\Yogu$
(or, more generally, $x\in N\cap G_0u$), 
the factorization $x=x_ux^u$ is given by
\[
\left[
\begin{array}{ccc}
1 & x_{12} & x_{13} \\
0 & 1 & x_{23} \\
0 & 0 & 1
\end{array}
\right]
=\left[
\begin{array}{ccc}
1 & x_{12} & 0 \\
0 & 1 & 0 \\
0 & 0 & 1
\end{array}
\right]
\left[
\begin{array}{ccc}
1 & 0 & x_{13}-x_{12}x_{23} \\
0 & 1 & x_{23} \\
0 & 0 & 1
\end{array}
\right]\,.
\]
The fiber of the projection $\pi_u:x\mapsto x_u$ over a point 
$x_u={\left[
\begin{array}{ccc}
1 & a & 0 \\
0 & 1 & 0 \\
0 & 0 & 1
\end{array}
\right]}\in \Yo_u$, $a>0$,
is therefore 
\[
\pi_u^{-1}(x_u)=
Y\cap\{ x_{12}=a\}
=
\left\{
\left[
\begin{array}{ccc}
1 & a & x_{13} \\
0 & 1 & x_{23} \\
0 & 0 & 1
\end{array}
\right]
\,:\,
ax_{23}\geq x_{13}\geq 0
\right\}\,. 
\]

\section{Transversals and links.
Proof of Theorem~\ref{th:strata}}
\label{sec:links}

Our next goal is to prove that the restriction of the projection 
$\pi_u$ onto $Y_{[u,v]}$ is globally trivialized along
$\Yo_u\,$.

\begin{figure}[ht]
\setlength{\unitlength}{1.8pt} 

\begin{center}
\begin{picture}(85,55)(-5,-7)

\qbezier(40,0)(30,30)(10,40)
\qbezier(40,0)(55,15)(60,40)

\thicklines
\qbezier(0,-5)(40,5)(80,-5)

  \put(5,27){\circle*{1.5}}
  \put(0,27){\makebox(0,0){$\tilde x$}}
\thinlines
\qbezier[50](5,27)(20,32)(35,30)

  \put(35,30){\circle*{1.5}}
  \put(40,0){\circle*{1.5}}

  \put(85,0){\makebox(0,-5){$\Yo_u$}}
  \put(85,20){\makebox(0,0){$\Yogu$}}
  \put(40,-5){\makebox(0,0){$x_u$}}
  \put(40,43){\makebox(0,0){$\pi_u^{-1}(x_u)$}}

  \put(40,26){\makebox(0,0){$\rho_{x_u}(\tilde x)$}}

\end{picture}
\end{center}

\caption{The projection $\rho_{x_u}$}
\label{fig:projection}
\end{figure}
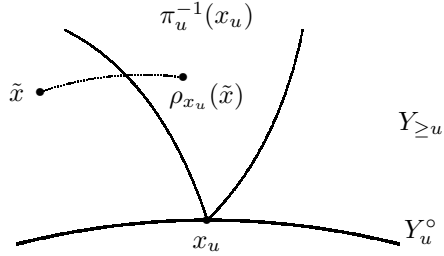

\begin{lemma}
\label{lem:hor-projection}
For any $\tilde x\in N\cap G_0 u$ and any $x_u\in N^u$, 
there exists unique $n_-\in N_-(u)$ such that 
the element $x'=[\tilde xn_-]_+$ is well defined and belongs to $x_u N(u)$. 

If moreover $\tilde x$ and $x_u$ are totally nonnegative,
then $x'$ is also totally nonnegative.
We thus obtain a cell-preserving projection 
\begin{eqnarray}
\label{eq:rho}
\begin{array}{rccc}
\rho_{x_u}:& \Yogu&\to&\pi_u^{-1}(x_u)\\[.1in]
           & \tilde x&\mapsto&x'  
\end{array}
\end{eqnarray}
(see Figure~\ref{fig:projection}). 
\end{lemma}

\begin{proof}
Let $\tilde x=\tilde x_u\tilde x^u$, where $\tilde x_u\in N^u$
and $\tilde x^u\in N(u)$, as in
Lemma~\ref{lem:pi_u}. 
Let $n_1\in N_-(u)$ be such that
$x_u=[\tilde x_u n_1]_+$ (such $n_1$ exists and is unique by
Lemma~\ref{lem:recovering-shift}). Set 
\begin{equation}
\label{eq:n_-}
n_-=[(\tilde x^u)^{-1} n_1]_- \,; 
\end{equation}
since both $\tilde x^u$ and $n_1$ belong to $u^{-1}Nu$, 
the element $n_-$ is well defined in view of Lemma~\ref{lem:wBw-in-G0}, 
and belongs to $N_-(u)$ by Lemma~\ref{lem:[]+insubgroup}.
Let us prove that the element $n_-$ defined by (\ref{eq:n_-})
has the desired properties, i.e., $x'=[\tilde xn_-]_+$
is well defined and belongs to $x_u N(u)$,
as shown below:  
\begin{equation}
\label{eq:comm-diagram}
\begin{array}{rcl}
\tilde x\,\, & \leadsto & x'=[\tilde xn_-]_+ \\[.1in]
\pi_u\,\downarrow\,\, & & \downarrow \\[.1in]
\tilde x_u & \leadsto & x_u=[\tilde x_u n_1]_+\,. 
\end{array}
\end{equation}
Denote $z=\tilde x^u n_-\,$.
Once again, $z\in u^{-1}Nu\subset N_-N\,$, and 
(\ref{eq:n_-}) implies 
$[z]_-=[\tilde x^u n_-]_-=[\tilde x^u [(\tilde x^u)^{-1} n_1]_-]_-
=n_1\,$. 
Then $\tilde x n_-=\tilde x_u z=\tilde x_u n_1 [z]_+\in G_0$
(because $[\tilde x_u n_1]_+=x_u$), so $x'$ is well defined indeed. 
Furthermore, 
$x'=x_u [z]_+\in x_u N(u)$
(by Lemma~\ref{lem:[]+insubgroup}), as desired. 

Uniqueness is proved by a similar argument. 
Suppose that $n_-\in N_-(u)$ is such that $x'=[\tilde xn_-]_+\in x_uN(u)$. 
As before, denote $z=\tilde x^u n_-\,$. 
Then $z=[z]_-[z]_+$ and $x'=[\tilde x_u z]_+
=[\tilde x_u[z]_-]_+\cdot [z]_+$. 
Since $[\tilde x_u[z]_-]_+\in N^u$ (by Lemma~\ref{lem:action}) 
and $[z]_+\in N(u)$, it follows from Lemma~\ref{lem:pi_u}
that $[\tilde x_u[z]_-]_+=x_u=[\tilde x_u n_1]_+\,$.
Hence by Lemmas~\ref{lem:[]+insubgroup}
and~\ref{lem:recovering-shift},  
$n_1=[z]_-=[\tilde x^u n_-]_-\,$, implying (\ref{eq:n_-}). 

It remains to prove the second part of the lemma. 
In view of Lemma~\ref{lem:recovering-shift}, 
a path connecting $\tilde x_u$ and $x_u$ within $\Yo_u$
gives a continuous deformation of the identity $1\in G$ into $n_1$
within $N_-(u)$,
which gives rise (via (\ref{eq:n_-}))
to a continuous deformation of $1$ into $n_-$ and, 
finally, to a path connecting 
$\tilde x$ and $x'=[\tilde xn_-]_+$ within the Bruhat cell 
containing~$\tilde x$. Hence $x'$ is totally nonnegative
by Proposition~\ref{prop:connected component}. 
\end{proof}

\subsection*{Example: $G=SL(3,\RR)$, $u=s_1\,$} 
For 
\[
\tilde x=\left[
\begin{array}{ccc}
1 & \tilde x_{12} & \tilde x_{13} \\
0 & 1 & \tilde x_{23} \\
0 & 0 & 1
\end{array}
\right]
\in\Yogu
\quad
\textrm{and} 
\quad
x_u=\left[
\begin{array}{ccc}
1 & a & 0 \\
0 & 1 & 0 \\
0 & 0 & 1
\end{array}
\right]\in \Yo_u\ , 
\]
computations give 
\[
n_-=
\left[
\begin{array}{ccc}
1 & 0 & 0 \\
a^{-1}-\tilde x_{12}^{-1} & 1 & 0 \\
0 & 0 & 1
\end{array}
\right]
\]
and
\[
x'=\rho_{x_u}(\tilde x)=[\tilde x n_-]_+
=\left[
\begin{array}{ccc}
1 & a & \displaystyle\frac{a\tilde x_{13}}{\tilde x_{12}} \\[.15in]
0 & 1 & \displaystyle\frac{\tilde x_{12}\tilde x_{23}-\tilde x_{13}}{a} 
+\displaystyle\frac{\tilde x_{13}}{\tilde x_{12}}\\[.15in]
0 & 0 & 1
\end{array}
\right]
\,. 
\]
Total nonnegativity of $x'$ does indeed follow from total
nonnegativity of $\tilde x$ and~$x_u\,$. 

\medskip

We denote by $\wnot$ the element of maximal length in~$W$. 

\begin{theorem}~\\ 
\label{th:transversal-sections}  
\begin{itemize}
\item[1.] 
For $x_u\in \Yo_u\,$, the set 
$x_u N(u)$ is a smooth submanifold in~$N\cap G_0 u$ 
diffeomorphic to the affine space $\RR^{\l(\wnot)-\l(u)}$. 
Furthermore, $x_u N(u)$ is 
transversal to every Bruhat stratum $N^w$, $w\geq u$ 
(hence every stratum $\Yo_w\subset\Yogu$). 

\item[2.]
For $x_u, \tilde x_u\in \Yo_u\,$, 
the map $\rho_{x_u}$ 
described in Lemma~\ref{lem:hor-projection} 
establishes a diffeomorphism 
between $\tilde x_u N(u)$ and $x_u N(u)$. 
This diffeomorphism respects total nonnegativity and 
the Bruhat stratification; more precisely, 
it restricts to a stratified diffeomorphism between the fibers 
$\pi_u^{-1}(\tilde x_u)$ and $\pi_u^{-1}(x_u)$. 
\end{itemize} 
\end{theorem}

\begin{proof}
The map $n_+\mapsto x_u n_+$ establishes a diffeomorphism between
$N(u)\cong \RR^{\l(\wnot)-\l(u)}$ and $x_u N(u)$. 
Let us prove transversality.
Consider a point $x\in x_u N(u)\cap N^w$. 
It will be enough to show that $x_u N(u)$ is transversal
to the smooth submanifold $[x N_-(u)]_+$ of dimension $\l(u)$ 
in $N^w$. 
Assume the contrary, i.e., there exists a common tangent vector
$v$ to $[x N_-(u)]_+$ and $x_u N(u)$ at the point~$x$. 
Let us evaluate the differential $D$ of the projection 
$N\cap G_0 u\to N^u$ at the vector~$v$. 
On the one hand, the projection is constant on $x_u N(u)$---hence 
$D(v)=0$. On the other hand, in view of (\ref{eq:n_-}),
the restriction of the projection
onto $[x N_-(u)]_+$ is a diffeomorphism---hence $D(v)\neq 0$, a
contradiction. 

Let us prove the second part of the theorem. From (\ref{eq:n_-})
and (\ref{eq:comm-diagram}), we have 
$x'=[\tilde x[(\tilde x^u)^{-1} n_1]_-]_+\,$,
where $\tilde x^u$ is given by (\ref{eq:x^u-formula})
and $n_1$ by Lemma~\ref{lem:recovering-shift}. 
The resulting map $\tilde x_u N(u)\to x_u N(u)$
is rational and therefore differentiable on its domain. 
Its inverse is again a map of the same kind,
with the roles of $x_u$ and $\tilde x_u$ reversed. 
Hence these maps are diffeomorphisms. 
Furthermore, they preserve the Bruhat stratification 
(in view of Lemma~\ref{lem:action}) 
and total nonnegativity (by the second part
of Lemma~\ref{lem:hor-projection}). 
\end{proof}

Recall the notation $Y_{[u,v]}=\bigcup_{w\in [u,v]} \Yo_w$
and $\Yogu=\bigcup_{w\geq u} \Yo_w$. 

\begin{corollary}
\label{cor:direct-product} 
For $u,v\!\in\! W$, $u\leq v$, 
and any $x_u\in \Yo_u\,$, 
we have the diffeomorphism of stratified spaces:
\[
Y_{[u,v]} \cong \Yo_u\times (\pi_u^{-1}(x_u))\cap Y_{[u,v]}) \,. 
\] 
In particular, $\Yogu\cong \Yo_u\times \pi_u^{-1}(x_u)$. 
\end{corollary}

\subsection*{Proof of Theorem~\ref{th:strata}}
Corollary~\ref{cor:direct-product} shows that the link 
of $\Yo_u$ in $Y_{[u,v]}$ is well defined  
(up to a stratified diffeomorphism),
and is explicitly given by
\[
\Lkuv = (\pi_u^{-1}(x_u))\cap Y_{[u,v]})\cap S_\varepsilon (x_u)\,,
\]
where $x_u$ is an arbitrary point on~$\Yo_u\,$, and 
$S_\varepsilon(x_u)$ is a small sphere centered at~$x_u\,$. 
The first two statements of Theorem~\ref{th:strata}
follow right away. 
The equalities (\ref{eq:closure+boundary}) 
follows from the analogous property
for the Bruhat stratification of~$Y$ (cf.\
Proposition~\ref{prop:lusztig-cells}), combined with  
Corollary~\ref{cor:direct-product}. 
\qed

\section{Proofs of 
Theorems~\ref{th:contractions} and~\ref{th:product}} 
\label{sec:proofs-type-A}


Recall that the elements $x_i(t)$ are defined by (\ref{eq:x_i(t)}). 
For the type~$A_{n-1}\,$, $x_i(t)$ is the $n\times n$
matrix that differs from the identity matrix in a single entry
(equal to~$t$) located in row $i$ and column~$i+1$.

\begin{definition} {\rm 
We define the regular map $\str:N\to\CC$ 
by the conditions
$\str(x_i(t))=t$ and $\str(xy)=\str(x)+\str(y)$. 
In particular, in the case of type~$A$, we have 
$\str(x)=\sum_i x_{i,i+1}\,$,
the sum of the matrix elements immediately above the main diagonal.
}\end{definition}

\begin{definition}{\rm 
For $\tau>0$, let $d(\tau)\in H$ be uniquely defined by 
the conditions $(d(\tau))^{\alpha_i}=\tau$, for all simple
roots~$\alpha_i\,$.
Then 
$d(\tau) x_i(a) d(\tau)^{-1}=x_i(\tau a)$ for any~$i$. 
For the type~$A_{n-1}\,$, 
\begin{equation}
\label{eq:d(t)}
d(\tau)=\tau^{-(n-1)/2}\left[
\begin{array}{cccc}
\tau^{n-1} & 0 & \cdots & 0 \\
0 & \tau^{n-2} & \cdots & 0 \\
\vdots & \vdots & \ddots & \vdots \\
0 & 0 & \cdots & 1
\end{array}
\right]\,, 
\end{equation}
and the automorphism $x\mapsto d(\tau)xd(\tau)^{-1}$ of the group $N$ 
multiplies each matrix entry $x_{ij}$ of~$x$ by~$\tau^{j-i}$. 
}\end{definition}

Note that the automorphism $x\mapsto d(\tau)xd(\tau)^{-1}$ 
preserves the cells $N^w$ and the subgroups $N(w)$;
it also preserves total nonnegativity.

\begin{definition}{\rm 
For $u\!\in\! W$ and $x_u\!\in\! \Yo_u\,$,
we define the vector field $\psi$ on $\pi_u^{-1}(x_u)$~by
\begin{equation}
\label{eq:psi}
\psi(x)=\frac{d}{d\tau}
(\rho_{x_u}(d(\tau)x d(\tau)^{-1}))\Bigl\vert_{\tau=1}
\end{equation}
(recall that $\rho_{x_u}$ is defined by (\ref{eq:rho})). 
}\end{definition}

\begin{lemma}
\label{lem:link}
{\rm (Type $A$ only.)} 
The vector field $\psi$ vanishes nowhere on $\pi_u^{-1}(x_u)$ except
at the point $x_u\,$. 
The directional derivative $\nabla_\psi\str(x)$ is positive
at every point $x\neq x_u\,$. 
\end{lemma}

The proof of this lemma is given in Section~\ref{sec:main-lemma}. 

\subsection*{Proof of Theorem~\ref{th:product}}
In view of Corollary~\ref{cor:direct-product}, 
it remains to show that the fiber $\pi_u^{-1}(x_u)\cap Y_{[u,v]}$ 
has the structure of the cone over the link $\Lkuv$. 

The vector field $\psi$ can be extended (by the same formula
(\ref{eq:psi}), with $x_u=\pi_u(x)$) 
to the open subset $N\cap G_0 u$ of~$N$. 
Furthermore, $\psi(x)$ is given by rational functions in the affine
coordinates of~$x$; therefore, the theorem of uniqueness and 
existence of solutions applies to this extension of~$\psi$
(hence to $\psi$ itself). 
Since $\psi$ is tangent to each stratum of the preimage 
$\pi^{-1}(x_u)$ (all these strata are smooth by Part~1 of 
Theorem~\ref{th:transversal-sections}),
it follows that every trajectory of $\psi$ is contained
in a single stratum.
 
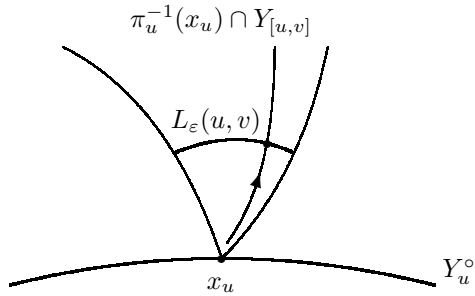
\begin{figure}[ht]
\setlength{\unitlength}{2pt} 

\begin{center}
\begin{picture}(85,55)(-5,-7)

\qbezier(40,0)(30,30)(10,40)
\qbezier(40,0)(55,15)(60,40)
\qbezier(41,3)(50,15)(50,40)

\thicklines
\put(47,15){\vector(1,3){0.5}}
\qbezier(0,-5)(40,5)(80,-5)
\qbezier(31,20)(43,25)(53.5,20)

  \put(40,0){\circle*{1.5}}
  \put(48.7,21.7){\circle*{1.5}}
  \put(39,26){\makebox(0,0){$L_\varepsilon(u,v)$}}

  \put(85,0){\makebox(0,-5){$\Yo_u$}}
  \put(40,-5){\makebox(0,0){$x_u$}}
  \put(40,45){\makebox(0,0){$\pi_u^{-1}(x_u)\cap Y_{[u,v]}$}}

\end{picture}
\end{center}

\caption{Embedding of the link into $\pi^{-1}(x_u)\cap Y_{[u,v]}$}
\label{fig:vector field}
\end{figure}

The intersection $Y\cap \{\str\leq c\}$ is compact for any $c>0$. 
Lemma~\ref{lem:link} then implies that for every $x_0\in \pi^{-1}(x_u)$, 
the solution $x(t)$ of the Cauchy problem $\dot x=\psi(x)$,
$x(0)=x_0$, $t<0$, exists for $t\in (-\infty,0]$. 
(Otherwise the trajectory $T_-=\{x(t):t\leq 0\}$ 
would hit the boundary of the stratum containing~$x_0\,$.) 
The trajectory $T_-$ must have limit points;
let $x_{\lim}$ be one of them. 
The function $s:t\mapsto \str(x(t))$, $t<0$,
is increasing (by Lemma~\ref{lem:link}) and bounded from below. 
Therefore $\lim_{t\to -\infty} \dot s(t)=0$,
implying $\nabla_\psi\str(x_{\lim})=0$. 
By Lemma~\ref{lem:link}, this means that $x_{\lim}=x_u\,$. 
Thus every trajectory of $\psi$ originates at the point $x_u$
(at $t=-\infty$).
A similar argument shows that $\lim_{t\to t^+} \str(x(t))=+\infty$,  
where $t^+$ denotes the upper limit of the maximal domain
of definition of $x(t)$ (so $t^+\in [0,\infty]$). 
We conclude that the function $\str$ increases from $\str(x_u)$ to
$\infty$ along each trajectory of~$\psi$,
except for the trajectory~$x(t)=x_u\,$. 
Thus every nontrivial trajectory $T\subset Y_{[u,v]}$ 
intersects the set 
\begin{equation}
\label{eq:L-epsilon}
L_\varepsilon(u,v)=L_{\varepsilon,x_u}(u,v)
=\pi_u^{-1}(x_u)\cap Y_{[u,v]} \cap 
\{ x: \str(x)=\str(x_u)+\varepsilon\} 
\end{equation}
at exactly one point; see Figure~\ref{fig:vector field}. 
Therefore $\pi^{-1}(x_u)\cap Y_{[u,v]}$ is diffeomorphic to the cone
${\rm Cone}(L_\varepsilon(u,v))$. 
(In particular, $\pi^{-1}(x_u)\cong {\rm Cone}(L_\varepsilon(u,\wnot).)$
This implies that $L_\varepsilon(u,v)$ is isomorphic
(as a stratified space) to the link~$\Lkuv$. 
Theorem~\ref{th:product} is proved.
\endproof

\subsection*{Proof of Theorem~\ref{th:contractions}}

For $x\in \pi_u^{-1}(x_u)\cap Y_{[u,v]}$, 
let $\lambda_{u,v}(x)$ denote the unique point of intersection of 
the link~$L_\varepsilon(u,v)$ with the trajectory of $\psi$ that
passes through~$x$
(see the sentence containing (\ref{eq:L-epsilon})). 

Fix $z\in Y_{\geq v}\,$, and 
define $R_{u,v}:L_\varepsilon(u,v)\times [0,1]\to L_\varepsilon(u,v)$ by
\begin{equation}
\label{eq:Ruv}
R_{u,v}(x,\tau)
=\lambda_{u,v}(\rho_{x_u}(\pi_v(d(\tau)zd(\tau)^{-1}
d(1-\tau)x d(1-\tau)^{-1})))\,.
\end{equation}
(Note that $d(\tau)zd(\tau)^{-1}\in Y_{\geq v}\,$,
and therefore $d(\tau)zd(\tau)^{-1}d(1-\tau)x d(1-\tau)^{-1}\in
Y_{\geq v}\,$, by Lemma~\ref{lem:YY-in-Y}. 
So the the right-hand side of (\ref{eq:Ruv}) is well defined.)

Let us show that the map $R_{u,v}$ is a deformation retraction 
of $L_\varepsilon(u,v)$ into a point. 
First, $R_{u,v}$ is continuous with respect to~$\tau$. 
Second, 
$\displaystyle\lim_{\tau\to 0} d(\tau)xd(\tau)^{-1}\!=\!1$ and
$\displaystyle\lim_{\tau\to 1} d(\tau)zd(\tau)^{-1}\!=\!z$.
Therefore $R_{u,v}(x,0)=x$.
On the other hand, 
$R_{u,v}(x,1)=\lambda_{u,v}(\rho_{x_u}(\pi_v(z)))$,
a point independent of~$x$.

It remains to observe that 
$\overline{S_{uvw}}\cong\overline{S_{uww}}\cong L_\varepsilon(u,w)$. 
\endproof

\section{Proof of Lemma~\ref{lem:link}}
\label{sec:main-lemma}

Our first goal is an explicit formula for the vector
field~$\psi$. 

Throughout this section, we use the following notation.
For $x\in \Yogu$, we denote 
\begin{eqnarray}
\label{eq:x_u,x^u,etc.}
\begin{array}{l}
x_u=\pi_u(x)=[u[u^{-1}[xu^{-1}]_+u]_-]_+\,,\\[.1in]
A=u^{-1} [x u^{-1}]_+u\,,\\[.1in]
y=[A]_-\,,\\[.1in]
x^u=[A]_+=x_u^{-1}x\,.\\[.1in]
\end{array}
\end{eqnarray}
Thus $x=x_u x^u$, $x_u=[uy]_+$
(agreeing with Lemma~\ref{lem:x and y}), and $A=y x^u$.

Let us fix an arbitrary totally nonnegative element $d\in H$. 
We will need some basic properties of the cell-preserving
automorphism $x\mapsto dxd^{-1}$. 

\begin{lemma}
\label{lem:pi_u(dxd)}
The automorphism $x\mapsto dxd^{-1}$ of $\Yogu$
commutes with the
maps $x\mapsto x_u$ and $x\mapsto x^u$:
\[
(dxd^{-1})_u=dx_ud^{-1}\,,\quad 
\left(dxd^{-1} \right)^u=dx^ud^{-1}\,.
\]
\end{lemma}

\begin{proof} 
In view of Lemma~\ref{lem:pi_u}, 
the statement follows from the factorization
$dxd^{-1}
=dx_ud^{-1}\cdot dx^ud^{-1}\,$, 
where the two factors on the right belong to $\Yo_u$ and
$N(u)$, respectively. 
\end{proof}

\begin{lemma}
\label{lem:y(dxd)}
The automorphism $x_u\mapsto dx_ud^{-1}$ of $\Yo_u$
commutes with the map 
$x_u\mapsto y$
(cf.\ Lemma~\ref{lem:x and y}).
In other words,
\[
dx_ud^{-1} = [u \cdot dyd^{-1}]_+ \,. 
\]
The unique element $n_1\in N_-(u)$
such that $x_u=[d x_u d^{-1} n_1]_+$ (cf.\
Lemma~\ref{lem:recovering-shift}) is given by $n_1=d y^{-1} d^{-1}y$. 
\end{lemma}

\begin{proof}
First part:
$d x_u d^{-1}=d [u y]_+ d^{-1}=
[d u y d^{-1}]_+=[ u d y d^{-1}]_+\,$. 
The second part is then a special case of
Lemma~\ref{lem:recovering-shift}. 
\end{proof}

\begin{lemma}
\label{lem:rho(xu)}
~$
\rho_{x_u}(d x d^{-1})
=
x \left(\left[ d A^{-1} d^{-1} A\right]_+\right)^{-1}
$. 
\end{lemma}

\proof 
Applying (\ref{eq:rho}) to $\tilde x=dxd^{-1}$
and using Lemmas~\ref{lem:pi_u(dxd)} and~\ref{lem:y(dxd)}
along with (\ref{eq:x_u,x^u,etc.}), 
we obtain: 
\[
\begin{array}{l}
\rho_{x_u}(d x d^{-1})
=
[ d x d^{-1}
[ d (x^u)^{-1} d^{-1}
d y^{-1} d^{-1}y]_-
]_+
=
[ d x d^{-1}
[ d A^{-1} d^{-1}A]_-
]_+ \\[.1in]
=
[ d x d^{-1}\cdot d A^{-1} d^{-1}A]_+
\left([ d A^{-1} d^{-1} A]_+\right)^{-1}
=
[x A^{-1} d^{-1}A]_+
\left([ d A^{-1} d^{-1} A]_+\right)^{-1}
\,.
\end{array}
\]
It remains to show that $[x A^{-1} d^{-1}A]_+=x$. 
This is done as follows:
\[
\begin{array}{l}
[x A^{-1} d^{-1}A]_+
=[x_u y^{-1}d^{-1}A]_+
=[[uy]_+y^{-1}d^{-1}A]_+
=[uy y^{-1}d^{-1}A]_+ \\[.1in]
=[ud^{-1}yx^u]_+
=[ud^{-1}y]_+ x^u
=[uy]_+ x^u
=x_u x^u
=x\,. \endproofmath
\end{array}
\]

Let ${\mathfrak g}$, ${\mathfrak n}$, and ${\mathfrak b_-}$ 
denote the Lie algebras of groups 
$G$, $N$, and $B_-\,$, respectively.
Let $\pi_{\mathfrak n}$ denote the projection 
${\mathfrak g}\to {\mathfrak n}$ along ${\mathfrak b}_-\,$. 
(In the case $\mathfrak{g}=\mathfrak{sl}_n$, 
$\pi_{\mathfrak n}$ replaces all lower-triangular entries
of a traceless matrix by zeroes.) 

\begin{lemma}
\label{lemma:derivative}
Let $f(\tau)=f_-(\tau)f_+(\tau)$, where $f_-(\tau)\in B_-\,$,
$f_+(\tau)\in N$ for all $\tau>0$.  
Assume that $f(1)=\operatorname{1}$. 
Then $f'_+(1)=\pi_{\mathfrak{n}}f'(1)$. 
\end{lemma}

\begin{proof}
The equality $f(1)\!=\!\operatorname{1}$ implies
$f_-(1)\!=\!f_+(1)\!=\!\operatorname{1}$. 
Then $f'(1)\!=\!f_-'(1)f_+(1)+f_-(1)f_+'(1)\!=\!f_-'(1)+f_+'(1)$.
Since $f_-'(1)\in {\mathfrak b_-}$ and 
$f_+'(1)\in {\mathfrak n}$, we are done. 
\end{proof}


\begin{proposition}
\label{prop:psi}
~$\psi(x) 
=x\pi_{\mathfrak n}\left(
A^{-1} d'(1) A\right)$. 
\end{proposition}

\proof 
By Lemma~\ref{lem:rho(xu)}, $\rho_{x_u}(d(\tau)xd(\tau)^{-1})
= x
\left([ d(\tau) A^{-1} d(\tau)^{-1}A]_+\right)^{-1}$. 
Hence 
\[
\psi(x)=\frac{d}{d\tau}
(\rho_{x_u}(d(\tau)xd(\tau)^{-1}))\Bigr\vert_{\tau=1}
=x
\frac{d}{d\tau}
\left(
[d(\tau) A^{-1} d(\tau)^{-1}A]_+
\right)^{-1}
\Bigr\vert_{\tau=1}\,.
\]
Applying Lemma~\ref{lemma:derivative}
and observing that $d(1)=1$, we obtain:
\[
\begin{array}{l}
\dfrac{d}{d\tau}
\left(
[d(\tau) A^{-1} d(\tau)^{-1}A]_+
\right)^{-1}
\Bigr\vert_{\tau=1}
=-
\frac{d}{d\tau}
[d(\tau) A^{-1} d(\tau)^{-1}A]_+ 
\Bigr\vert_{\tau=1}
\\[.2in]
=-\!
\pi_{\mathfrak n}
\frac{d}{d\tau}
\left( 
d(\tau) A^{-1} d(\tau)^{-1}A
\right)
\Bigr\vert_{\tau=1}
=-\pi_{\mathfrak n}
(d'(1) \!-\! A^{-1} d'(1) A)
=\pi_{\mathfrak n}
(A^{-1} d'(1) A),
\end{array}
\]
implying the claim.
\endproof

In the rest of this section, we only consider the case of the
type~$A_{n-1}\,$. 
Thus $W$ is the symmetric group~$\mathcal{S}_n\,$.
We treat the elements of $W$ as bijective maps
$\{1,\dots,n\}\to\{1,\dots,n\}$, 
and choose permutation matrices as their 
representatives in~$G=GL(n)$.



\begin{lemma} 
\label{lemma:pos_prod}
{\ }
\begin{enumerate} 
\item If $1\le i\le j\le n$
and $u(j)\leq u(i)\leq u(j+1)$, then 
$\left(A^{-1}\right)_{j,i}\cdot A_{i,j+1}\geq 0$.
\item If $1\le j<i\le n$
and $u(j)\leq u(i)\leq u(j+1)$, then
$\left(A^{-1}\right)_{j,i}\cdot A_{i,j+1}\leq 0$.
\item Otherwise, $\left(A^{-1}\right)_{j,i}\cdot A_{i,j+1}=0$.
\end{enumerate}
\end{lemma}

\begin{proof}
For a matrix $a\in u^{-1}N u$, the matrix element $a_{ij}$ vanishes
unless $u(i)\leq u(j)$. 
It follows that 
$\left(A^{-1}\right)_{j,i}\cdot A_{i,j+1}=0$
unless $u(j)\leq u(i)\leq u(j+1)$,
proving Part~3 of the lemma. 


Let us  prove Parts~1 and~2. 
Set $s=u(i)$, $p=u(j)$, $q=u(j+1)$; thus $p\leq s\leq q$.

In what follows, we denote by 
$z^{i_1,\dots,i_r}_{j_1,\dots,j_r}$ the determinant of the 
submatrix of a matrix $z$ formed by the rows 
$i_1,\dots,i_r$ 
and the columns $j_1,\dots,j_r$ (in this order). 
Using the definition of $A$ and the fact that $x\in N$, we obtain: 
\[
\begin{array}{l}
A_{i,j+1}
=([xu^{-1}]_+)_{s,q}
=([xu^{-1}]_+)^{1,\dots,s}_{1,\dots,s-1,q}\\[.1in]
=\displaystyle\frac{(xu^{-1})^{1,\dots,s}_{1,\dots,s-1,q}}{(xu^{-1})^{1,\dots,s}_{1,\dots,s}}
=\displaystyle\frac{x^{1,\dots,s}_{u^{-1}(1),\dots,u^{-1}(s-1),u^{-1}(q)}}{x^{1,\dots,s}_{u^{-1}(1),\dots,u^{-1}(s)}}\,. 
\end{array}
\]
Since $x$ is totally nonnegative, the sign of $A_{i,j+1}$ is either
zero or $(-1)^{n_1+n_2}$, 
where 
\[
\begin{array}{l}
n_1=\Card\{l\,:\,l\le s-1, u^{-1}(l)>i\}\,,\\[.1in]
n_2=\Card\{l\,:\,l\le s-1, u^{-1}(l)>j+1\}\,.
\end{array}
\]

The sign of $\left( A^{-1}\right)_{j,i}$ can be determined in a
similar fashion. 
Using the notation $\hat i$ to indicate that the index $i$ is being
removed, we obtain: 
\[
\begin{array}{l}
(A^{-1})_{j,i}
=(([xu^{-1}]_+)^{-1})_{p,s}
=(-1)^{(s-p)}
([xu^{-1}]_+)^{1,\dots,\hat s,\dots,n}_{1,\dots,\hat p,\dots,n}\\[.1in]
=(-1)^{(s-p)}
([xu^{-1}]_+)^{1,\dots,s-1}_{1,\dots,\hat p,\dots,s}
=(-1)^{(s-p)}\,\displaystyle\frac{(xu^{-1})^{1,\dots,s-1}_{1,\dots,\hat p,\dots,s}}{(xu^{-1})^{1,\dots,s-1}_{1,\dots,s-1}}\\[.1in]
=(-1)^{(s-p)}\,\displaystyle\frac{x^{1,\dots,s-1}_{u^{-1}(1),\dots,\widehat{u^{-1}(p)},\dots,u^{-1}(s)}}{x^{1,\dots,s-1}_{u^{-1}(1),\dots,u^{-1}(s-1)}}\,. 
\end{array}
\]
Hence the sign of $(A^{-1})_{j,i}$ is either
zero or $(-1)^{s-p+n_1+n_3+n_4}$,
where $n_1$ is the same as before, and 
\[
\begin{array}{l}
n_3=\Card\{l\,:\,l\le p-1, u^{-1}(l)>j\}\,,\\[.1in]
n_4=\Card\{l\,:\,p<l\le s, u^{-1}(l)<j\}\,.
\end{array}
\]
Since
\[
n_3+s-p-n_4=\Card\{l\,:\,l\le s, u^{-1}(l)>j\}\,,
\]
we conclude that the sign of $A_{i,j+1}(A^{-1})_{j,i}$ is either
zero or 
$(-1)^{\Card\{l\,:\,l=s, u^{-1}(l)>j\}}$,
matching the claim of the lemma. 
\end{proof}

Let us extend the notation $\str(x)=\sum_i x_{i,i+1}$ 
to arbitrary matrices~$x$.

\begin{lemma}
\label{lemma:str}
Let $\nu=\emph{\texttt{diag}}(n,n-1,\dots,1)$. 
Then $\emph{\str}( A^{-1}\nu\,A)\geq 0$.
\end{lemma}

\begin{proof} 
For $k=1,\dots,n$, 
let $\nu_k=\texttt{diag}(1,\dots,1,0,\dots,0)$ 
denote the diagonal matrix
whose first $k$ diagonal entries are equal to~1, 
and all other entries vanish.
The equality $A^{-1}A=1$ implies 
\[
\sum_{i=1}^j 
\left( A^{-1}\right)_{j,i}\cdot A_{i,j+1}+
\sum_{i=j+1}^n 
\left( A^{-1}\right)_{j,i}\cdot A_{i,j+1}=0\,,
\]
where by Lemma~\ref{lemma:pos_prod} all terms in the first sum
are nonnegative while all terms in the second sum are nonpositive.
Then
$\left( A^{-1}\nu_k A\right)_{j,j+1}=\sum_{i=1}^k 
\left( A^{-1}\right)_{j,i}\cdot A_{i,j+1}\geq 0\,$,  
implying $\str(A^{-1}\nu_k A)\ge 0$. 
Since $\nu=\sum_{k=1}^n \nu_k\,$, the lemma follows.
\end{proof}

We are now prepared to complete the proof of Lemma~\ref{lem:link}.
First, let us note that $d(\tau)x_u d(\tau)^{-1}\in \Yo_u\,$;
hence $\rho_{x_u}(d(\tau)x_u d(\tau)^{-1})\in \Yo_u\cap
\pi_u^{-1}(x_u)=\{x_u\}\,$. 
Thus $\psi(x_u)=0$.

Let $x\in \pi^{-1}(x_u)$, $x\neq x_u\,$. 
Since $d'(1)=\nu+\lambda$, for some scalar matrix~$\lambda$, 
Proposition~\ref{prop:psi} yields  
$\psi(x) 
=x\pi_{\mathfrak n}\left(
A^{-1} \nu A\right)$. 
Therefore 
\[
\str(\psi(x))
=\str(x)+\str(\pi_{\mathfrak n}(A^{-1} \nu A))
=\str(x)+\str(A^{-1} \nu A)
\]
(here we identify the tangent vector $\psi(x)$ with the corresponding
traceless matrix). 
Since $\str(x)>0$, Lemma~\ref{lemma:str} implies that
$\str(\psi(x))>0$; in particular, $\psi(x)\neq 0$. 
Then $\str(x+\psi(x)d\tau)=\str(x)+\str(\psi(x))d\tau$;
hence $\nabla_\psi \str(x)=\str(\psi(x)) >0$, as desired.
\endproof

\end{document}